\documentclass[a4paper,10pt]{article}
\usepackage{amsmath}
\usepackage{amsthm}
\usepackage{hyperref}
\usepackage{amsfonts}
\usepackage{amssymb}
\usepackage{latexsym}
\newtheorem{theorem}{Theorem}[section]
\newtheorem{lemma}[theorem]{Lemma}

\title{\large{\textbf{ON THE ASYMPTOTIC EXPANSION OF THE \\ SUM OF THE FIRST $n$ PRIMES}}}

\author{\normalsize{\textit{Nilotpal Kanti Sinha}}}
 
 \date{}
 
\begin{document}

\maketitle
{
\begin{flushleft}
\textit{Dedicated to Dr. Madan Mohan Singh, my first mathematics teacher.}
\end{flushleft}
}

\begin{abstract}
\small{An asymptotic formula for the sum of the first $n$ primes is derived. This result improves the previous results of P. Dusart. Using this asymptotic expansion, we prove the conjectures of R. Mandl and G. Robin on the upper and the lower bound of the sum of the first $n$ primes respectively.}
\end{abstract}

\section{Introduction}

Let $p_n$ denote the $n^{th}$ prime 
\footnote{    2000  \textit{Mathematics Subject Classification.} 11A41.}
\footnote{  \textit{Key words and phrases.} Primes, Inequalities.}.
Robert Mandl conjectured that  
\begin{displaymath}
\sum_{r \le n}p_r < \frac{np_n}{2}.
\label{mandl}
\end{displaymath}

This conjecture was proven by Rosser and Schoenfeld in \cite{rs} and is now referred to as Mandl's inequality. An alternate version of the proof was given by Dusart in \cite{pd}. In the same paper, Dusart also showed that
\begin{equation} 
\sum_{r \le n}p_r = \frac{n^2}{2}(\ln n + \ln\ln n - 3/2 + o(1)).
\label{dusart1}
\end{equation}

Currently, the best upper bound on Mandl's inequality is due to Hassani who showed that (See \cite{mh}) for for $n \ge 10$,  
\begin{equation}
\sum_{r \le n}p_r < \frac{np_n}{2} - \frac{n^2}{14}.
\label{hasanni}
\end{equation}

With regards to the lower bound, G. Robin conjectured that 
\begin{equation}
n p_{[n/2]} < \sum_{r <n}p_r.
\label{robin}
\end{equation}

This conjecture was also proved by Dusart in \cite{pd}. However, neither \ref{hasanni} nor \ref{robin} give exact growth rate of 
$\sum_{r \le n}p_r$. In this paper, we shall derive the asymptotic formula for $\sum_{r \le n}p_r$. Both Hassani's improvement of Mandl's inequality and Robin's conjecture follow as corollaries of our asymptotic formula.

\section{Asymptotic expansion of $\sum_{r \le n}p_r$}

\begin{theorem} (\textbf{M. Cipolla})
There exists a sequence $(P_m)_{m \ge 1}$ of polynomials with rational coefficients such that, 
for any integer m,
\begin{displaymath}
p_n = n \Bigg[\ln n + \ln\ln n - 1  + \sum_{r=1}^{m}\frac{(-1)^rP_r (\ln\ln n)}{\ln^r n} + 
o\Bigg(\frac{1}{\ln^m n}\Bigg)\Bigg].
\label{cipolla1}
\end{displaymath}
\end{theorem}

This was proved by M. Cipolla in a beautiful paper (See \cite{mc}) in 1902. In the same paper, Cipolla gives recurrence formula for $P_m$ and shows that every $P_m$ has degree $m$ and leading coefficient $\frac{1}{m}$. In particular,

\begin{equation}
P_1 (x) = x-2, P_2 (x) = \frac{1}{2}(x^2 - 6x + 11).
\label{cipolla2}
\end{equation}

\begin{lemma}
If f is monotonic and continuous and defined in $[1,n]$ and then,
\begin{displaymath}
\sum_{r\le n}f(r)= \int_{1}^{n} f(x)dx + O(|f(n)|+|f(1)|).
\label{lemma2.2}
\end{displaymath}
\end{lemma}

\begin{proof}
Well known. (See \cite{ik}, 1.62-1.67, Page 19-20)
\qedhere\
\end{proof}

\begin{theorem}
There exists a sequence $(S_m)_{m \ge 1}$ of polynomials with rational coefficients such that, 
for any integer m,
\begin{displaymath}
\sum_{r \le n}p_r = \frac{n^2}{2} \Bigg[\ln n + \ln\ln n - \frac{3}{2}  + \sum_{r=1}^{m}\frac{(-1)^{r+1}S_r (\ln\ln n)}{r\ln^r n}
+ o\Bigg(\frac{1}{\ln^m n}\Bigg)\Bigg].
\end{displaymath}
Further, every $S_m$ has degree m and leading coefficient $1/m$. In particular
\\
\begin{displaymath}
S_1 (x) = x-\frac{5}{2}, S_2 (x) = x^2 - 7x + \frac{29}{2}.
\end{displaymath}
\end{theorem}

\begin{proof} 
We define $p(x)$ as
\begin{equation}
p(x) = x\ln x + x\ln\ln x - x  + \sum_{r=1}^{m}\frac{(-1)^rxP_r (\ln\ln n)}{\ln^r n}.
\label{theorem2.3a}
\end{equation}
where $P_r(x)$ is the same sequence of polynomials as in Theorem \ref{cipolla1}. It follows form Lemma \ref{lemma2.2} that

\begin{displaymath}
\sum_{r \le n}p_r = 2 + 3 + \int_{3}^{n}p(x)dx + O(p_n) + o\Bigg(\int_{3}^{n}\frac{x}{\ln^m x}\Bigg).
\end{displaymath}
Now $p_n \sim n \ln n$ where as

\begin{displaymath}
\int_{3}^{n}\frac{x}{\ln^m x} \sim \frac{n^2}{2\ln^m n}
\end{displaymath}
which grows much faster than $n \ln n$. Hence

\begin{equation}
\sum_{r \le n}p_r = \int_{3}^{n}p(x)dx + o\Bigg(\frac{n^2}{\ln^2 n}\Bigg).
\label{theorem2.3b}
\end{equation}

All we need to do is the integrate each term of \ref{theorem2.3a}. Except for a couple of simple terms, integration of the terms of \ref{theorem2.3a} will result in an infinite series and due to \ref{theorem2.3b}, we can stop the series when the growth rate of a new term is equal to or slower than the error term in \ref{theorem2.3b}. Since $P_m(\ln\ln x)$ is a polynomial of degree $m$ and has rational coefficients with leading coefficient $1/m$, the integration of each terms of $p(x)$ will result in an infinite series of the type

\begin{displaymath}
\int_{3}^{n}xP_m(\ln\ln x)dx = \frac{(-1)^m n^2}{2}\sum_{i=1}^{\infty}Q_{m,i}(\ln\ln n) + O(1)
\end{displaymath}
where $Q_{m,i}(x)$ is a polynomial of degree $m$ with rational coefficients and leading coefficient $1/m$. Thus the polynomial $S_m(x)$ is of degree $m$ and has rational coefficients with leading coefficient $1/m$.
\\

To find the first two terms of the polynomial $S_m(x)$ we integrate the first four terms of $p(x)$. The first four terms of $p(x)$ are

\begin{equation}
x\ln x + x\ln\ln x - x  + \frac{x\ln\ln x - 2x}{\ln x} - \frac{x\ln^2 \ln x -6x\ln\ln x + 11x}{2\ln^2 x}.
\label{theorem2.3c}
\end{equation}
Integrating each term separately, we have
 
\begin{equation}
\int_{3}^{n} x \ln x dx = \frac{n^2 \ln n}{2} - \frac{n^2}{4} + O(1)
\label{theorem2.3d}
\end{equation}

\begin{equation}
\int_{3}^{n} x \ln\ln x dx = \frac{n^2 \ln\ln n}{2} - \frac{n^2}{4\ln n} - \frac{n^2}{8\ln^2 n}
+ O\Bigg(\frac{n^2}{\ln^3 n}\Bigg)
\end{equation}

\begin{equation}
-\int_{3}^{n} x  dx = -\frac{n^2}{2} + O(1)
\end{equation}

\begin{equation}
\int_{3}^{n} \frac{x\ln \ln x}{\ln x} dx = \frac{n^2 \ln \ln n}{2 \ln n} + \frac{n^2 \ln \ln n}{4 \ln^2 n}
- \frac{n^2}{4 \ln n} + O\Bigg(\frac{n^2 \ln \ln n}{\ln^3 n}\Bigg)
\end{equation}

\begin{equation}
- 2\int_{3}^{n} \frac{x}{\ln x}dx = - \frac{n^2}{\ln n} - \frac{n^2}{2\ln^2 n} + O\Bigg(\frac{n^2}{\ln^3 n}\Bigg)
\end{equation}

\begin{equation}
-\frac{1}{2}\int_{3}^{n} \frac{x\ln^2 \ln x}{\ln^2 x} dx = -\frac{n^2 \ln^2 \ln n}{4 \ln^2 n} 
+ O\Bigg(\frac{n^2 \ln^2 \ln n}{\ln^3 n}\Bigg)
\end{equation}

\begin{equation}
3 \int_{3}^{n} \frac{x\ln\ln x}{\ln^2 x} dx = \frac{3n^2 \ln\ln n}{2\ln^2 n} 
+ O\Bigg(\frac{n^2 \ln^2 \ln n}{\ln^3 n}\Bigg)
\end{equation}

\begin{equation}
-\frac{11}{2} \int_{3}^{n} \frac{x}{\ln^2 x} dx = -\frac{11n^2}{4\ln^2 n} 
+ O\Bigg(\frac{n^2 \ln n}{\ln^3 n}\Bigg)
\label{theorem2.3e}
\end{equation}
Adding \ref{theorem2.3d} - \ref{theorem2.3e} we obtain

\begin{displaymath}
\sum_{r \le n}p_r = \frac{n^2}{2}\Bigg[\ln n + \ln\ln n - \frac{3}{2} + \frac{\ln\ln n}{\ln n} - \frac{5}{2\ln n}
- \frac{\ln^2 \ln n}{2\ln^2 n}
\end{displaymath}

\begin{equation}
+ \frac{7 \ln \ln n}{2\ln^2 n} - \frac{29}{4\ln^2 n}
+ o\Bigg(\frac{1}{\ln^2 n}\Bigg) \Bigg].
\label{theorem2.3f}
\end{equation}
\\
Notice that taking the first four terms of \ref{theorem2.3f} we obtain Dusart's result in \ref{dusart1}. This proves the theorem.
\qedhere\
\end{proof}

\section{The inequality of Robin}

From the asymptotic expansion of $\sum_{r \le n}p_r$ we can not only prove the inequalities of Mandl \ref{mandl}and Robin \ref{robin} but also refine them.

\begin{lemma}
\begin{displaymath}
\sum_{r < n}p_r = n p_{[n/2]} + \frac{2\ln 2 - 1}{4}n^2 + O\Bigg(\frac{n^2 \ln\ln n}{\ln n}\Bigg).
\label{lemma3.2}
\end{displaymath}
\end{lemma}

\begin{proof}
Taking $[n/2]$ in place of $n$ in the asymptotic expansion of $n^{th}$ prime we obtain
\begin{displaymath}
n p_{[n/2]} = \frac{n^2}{2}(\ln n + \ln \ln n - 1 - \ln 2) + O\Bigg(\frac{n^2 \ln\ln n}{\ln n}\Bigg)
\end{displaymath}
\begin{displaymath}
 = \frac{n p_n}{2} -\frac{n^2 \ln 2}{2} + O\Bigg(\frac{n^2 \ln\ln n}{\ln n}\Bigg)
\end{displaymath}
Using \ref{theorem2.3f} we can reduce this to
\begin{displaymath}
 = \sum_{r <n}p_r + \frac{n^2}{4} -\frac{n^2 \ln 2}{2} + O\Bigg(\frac{n^2 \ln\ln n}{\ln n}\Bigg)
\end{displaymath}
This proves the lemma.
\qedhere\
\end{proof}

Since the second term of Lemma \ref{lemma3.2} is positive, it follows that for all sufficiently large $n$, Robin's conjecture is true.

\small{e-mail: 
\texttt{nilotpalsinha@gmail.com, nilotpal.sinha@greatlakes.edu.in}}
\end{document}